\documentclass[english]{llncs}
\usepackage{lmodern}

\usepackage[T1]{fontenc}
\usepackage[latin9]{inputenc}
\usepackage{babel}
\usepackage{array}
\usepackage{varioref}
\usepackage{float}
\usepackage{multirow}
\usepackage{graphicx}
\usepackage[unicode=true,
 bookmarks=true,bookmarksnumbered=false,bookmarksopen=false,
 breaklinks=true,pdfborder={0 0 0},backref=false,colorlinks=false]
 {hyperref}
\hypersetup{pdftitle={Set Theory for The (Smart) Masses},
 pdfauthor={Moez A. AbdelGawad},
 pdfsubject={Set Theory in Android Proof Maker},
 pdfkeywords={Set Theory; Proof Designer; Smartphones; Proof Maker; Proof Assistants; Coq; Isabelle}}

\makeatletter

\providecommand{\tabularnewline}{\\}

\newcommand{\code}[1]{\texttt{#1}}

\makeatother

\begin{document}

\title{Set Theory for The (Smart) Masses}

\author{Moez A. AbdelGawad\\
moez@cs.rice.edu}

\institute{SRTA-City, Alexandria, Egypt}
\maketitle
\begin{abstract}
Proof Designer is a computer software program designed to help Mathematics
students learn to write mathematical proofs. Under the guidance of
the user, Proof Designer assists in writing outlines of proofs in
elementary set theory. Proof Designer was designed by Daniel Velleman
in association with his book \textquotedbl{}How To Prove It: A Structured
Approach\textquotedbl{} to help students apply the methods discussed
in the book, making classes based on the book more interactive.

This paper is an early report on the progress of our effort to ``bring
set theory to the masses'' by developing Proof Maker, a new Proof
Designer-inspired software that ports Proof Designer to hand-held
devices such as smart-phones and tablets. Proof Maker, when completed,
will allow students to use Proof Designer with the ease of a touch,
literally, on their smart devices. Our goal behind developing Proof
Maker is to enable any one who is interested enough to develop elementary
set theory proofs anywhere he or she might be (think of doing proofs
while waiting at a bus stop!) and at any time he or she wishes (think
of writing proofs before going to bed, or even in bed!). In this paper
we report on the improvements we made to Proof Designer so far, and
on the (many) steps remaining for us to have a fully-functioning Proof
Maker ``in our hands''.
\end{abstract}

\section{Introduction}

Mathematics is useful for all branches of scientific research. The
mastery of mathematical skills is an essential enabler of success
in almost all sciences. The mastery of discrete mathematics, which
studies discrete and distinct mathematical objects, is particularly
important for many branches of scientific research, including, for
example, the efficient production of correct and efficient software.

Unfortunately, globally-speaking, the mathematical skills of undergraduate
and even graduate students are significantly lacking. Many prospective
scientists lack the basics of how to think mathematically and how
to write a correct mathematical proof, despite the importance of such
skills for their future success as scientists.

To help in addressing this situation, Daniel Velleman wrote a book
in 1994 (with a second edition in 2006) titled `How To Prove It: A
Structured Approach'~\cite{Velleman06}, in which he likened constructing
mathematical proofs to structured programming. Velleman used in his
book examples from arithmetic and high-school mathematics to present
his ideas on how to construct mathematical proofs in a structured
way.

Velleman encouraged the use of his book by referring to a pedagogic
Mathematics software that he developed called Proof Designer, by which
readers of his book can apply the ideas they learn from the book.
Proof Designer is freely-available online, as a Java applet, to assist
its users build mathematical proofs of elementary set theory theorems
in a structured way. Since 1994, and more so since 2006, many mathematics
courses around the world have used Velleman's book, and its accompanying
software, as an essential references for teaching the skills of mathematical
thinking and proof construction to graduate and undergraduate students.
Helping in the widespread use of the book was Velleman's lucid writing
style, his use of elementary mathematical examples in his book, and
also the ease of use of Proof Designer when compared to that of other
proof assistants.

Recently we started a project whose goal is to take Proof Designer
to its next step, so as to make it usable in wider contexts and to
appeal to an even wider audience. Nowadays, in the age of handheld
devices (such as tablets and smart-phones), Proof Designer is starting
to show its age and limitations. For example, there is no portal of
Proof Designer to any of the popular platforms for handheld devices.
Additionally, Proof Designer uses only English as the language of
its proofs and the language of its graphical user interface. Thus,
compared to GUIs of modern educational software, despite its success
and it fully serving its initial purpose, Proof Designer is now clearly
lacking in many regards. The goal of our project is to eliminate most,
if not all, of the limitations on Proof Designer that make it less-used
today as a math education software than it was during the last ten
years. In the following two sections we report, using software illustrations,
on our effort so far. We first describe in the next section what we
have done so far, then, in the following section, we describe what
remains to be done.

\section{What Has Been Done}

To describe what we have done so far, we first describe Proof Designer
in its original form then describe changes we made to it. Then we
describe steps we made so far towards porting Proof Designer to the
Android platform.

\subsection{Original Proof Designer}

Proof Designer allows users to develop proofs using an intuitive interface.
Figures~\ref{fig:An-Incomplete-Proof}-\ref{fig:Proof-Designer-Reexpress},
on pp.~\pageref{fig:An-Incomplete-Proof}-\pageref{fig:Proof-Designer-Reexpress},
show the main components of the Proof Designer user experience, which
involve the presentation of structured complete and incomplete proofs,
a theorem-entry dialog box, drop-down menus that the user uses to
construct his or her proofs, and a dialog box for re-expressing mathematical
formulas.

Proof Designer, in its original form, is available for use as a Java
applet at http://www.cs.amherst.edu/\textasciitilde{}djv/pd/pd.html.
Instructions for how to setup and use Proof Designer can be found
at http://www.cs.amherst.edu/\textasciitilde{}djv/pd/help/Instructions.html.

\subsection{Proof Designer Improvements}

Before setting on building Proof Maker as a portal of Proof Designer
to handheld devices, we set on making some improvements to Proof Designer
itself. After communicating with Professor Velleman and consulting
with him, he kindly sent us the source code of Proof Designer. We
made many changes to Proof Designer, some of which are visible to
the user, and some are not.

\subsubsection{Code Improvements (Invisible to User)}

First, to enhance our understanding of the Proof Designer code base
and to facilitate its further development, we made some improvements
to the software source code. In particular,
\begin{enumerate}
\item Proof Designer had little documentation for its source code. We thus
added some unit tests, \code{assert} statements, and code comments.
\item \noindent All classes of Proof Designer were in one Java package (the
default package). Based on UML class diagrams of the code base (see
Figure~\vref{fig:Class-hierarchy-PD}), we distributed the code among
seven Java packages (a.k.a., ``modules''), the most important being
packages for formula classes (class \code{Formula} and its descendant
classes) and for proof component classes (class \code{PComponent}
and its descendants). We also had packages for class \code{MenuAction}
and all its \code{DoX} descendant classes, and for class \code{PDialog}
and all its descendant Proof Designer dialog box classes (e.g., class
\code{EntryDlg}, which is used to enter theorem statements in Proof
Designer).
\item Proof Designer was written using Java 1.3. Hence, its code made no
use of Java generics or Java enumerations, for example, which were
introduced in Java 1.5/5.0. We thus used generics wherever possible
in the Proof Designer source code to improve the reliability and maintainability
of the code, and we also made use of \code{Enum}s (instead of \code{int}s),
e.g., for defining Proof Designer's formula and operator kinds.
\end{enumerate}

\subsubsection{Visible Improvements}

We also made changes that are visible to the user, to improve his
or her user experience.
\begin{enumerate}
\item \noindent We restructured menus so that some user actions, more intuitively,
are viewed as either inferences (from givens) or are goal-oriented
actions. (See Figure~\vref{fig:NPD-Menus}.)
\item We added the ability to save and load proof sessions (as XML files).
\item We added the ability to run Proof Designer, not only as a web browser
applet but also as a standalone software (a Java jar file) that can
be downloaded and run without the need for a web browser.
\item Proof Designer originally had the ability to a single undo/redo proof
step. We added an unlimited undo/redo capability to Proof Designer.
\item \noindent We also added a new unlimited undo/redo capability in Proof
Designer's Reexpress dialog. (See Figure~\vref{fig:NPD-Reexpress}.)
\item We added limited support for automating proofs in Proof Designer by
adding an `Auto' command for use on proof goals. The Auto command
automatically decides and performs the next step in the proof, if
any, based on the logical form of the goal statement.
\item To ease the use of Proof Designer (and to gear it more towards touch-based
interaction), we added a toolbar that has the auto and undo/redo commands.
(See Figure~\vref{fig:NPD-Toolbar}.)
\end{enumerate}

\subsection{APM (Android Proof Maker): Porting Proof Designer to Android}

After implementing the above-mentioned improvements to Proof Designer,
we set on exploring porting Proof Designer to handheld devices. We
decided to call the new software Proof Maker. Given the global widespread
use of the Android platform, we picked the platform as our first choice
for porting Proof Designer to. We call the portal to the Android platform
Android Proof Maker (or, APM for short).

Given that typical Android software is written using Java, we initially
assumed porting Proof Designer to Android will be straightforward.
In fact the \code{Formula} package in Proof Designer (after making
the above-mentioned changes) was ported without a single change to
its code. However, we soon realized that there is no one-to-one correspondence
between Java Swing UI (user interface) components (used in Proof Designer)
and Android UI components. The differences include, for example,
\begin{itemize}
\item The Android \code{View} class has a somewhat different semantics
and a different behavior than the Java Swing \code{JComponent} class.
\item The Android's \code{ViewGroup} class is different from its Swing
approximate counterpart class \code{Container}.
\item Although the Android platform has dialog boxes, but the closest to
a Java Swing dialog box is usually an Android \code{Activity} not
an Android dialog box.
\item Similarly, \code{JFrame} and \code{JPanel} in Java Swing have no
exact counterparts in Android UI components. The closest Android classes
to them seem to be \code{Activity} and \code{LinearLayout}, respectively.
\item In Java Swing the \code{Toolkit} and \code{Font} classes provide
font services that in Android are provided, using a different API,
in classes \code{Paint} and \code{TypeFace}.
\end{itemize}
We thus started experimenting with Android UI components to see which
could suit our purposes and best approximate the Proof Designer user
interface, and that will incur the least changes to the source code
of Proof Designer so as to maintain as much as possible of its ``spirit''.
(See Figure~\vref{fig:APM-Theorem-Entry} and Figure~\vref{fig:APM-Incomplete-Proof}.)
Even though not as polished as their Proof Designer counterparts,
our portal of some of the main Proof Designer UI components to the
Android platform is a good proof of concept that the portal is possible,
even when it will not be straightforward. Our effort so far has provided
us thus with an assurance that Proof Designer does not need a total
rewriting to be ported to the Android platform or to platforms of
other handheld devices.

It is worthy to mention that due to Proof Designer not employing the
popular Model-View-Controller (MVC) model in its software design,
we do though expect the differences between the Android and Swing
UI APIs to affect the final versions of Android Proof Maker, particularly
affecting the presentation of proofs, which in Proof Designer are
modeled using descendants of class \code{PComponent}. (Contrary to
the requirements of MVC, class \code{PComponent} and its descendants
in Proof Designer doubly function as Model classes, modeling abstract
proof components, but also as View classes that inherit from the Swing
\code{JComponent} UI class and as such are used as part of the GUI
of Proof Designer).

\section{What Remains to Be Done}

As demonstrated by the figures for APM as we have it today, Proof
Maker is still far from complete. Much remains to be done before we
get to a final usable version of Proof Maker. We mention the most
important remaining steps below.

\subsection{To Be Done in Proof Designer}

We first intend to make further improvements to Proof Designer. These
include the following.

\subsubsection{Code Improvements}
\begin{enumerate}
\item Adding more unit tests, assertions, and code comments.
\item Consider using the MVC software design model. Mostly will affect proof
components (\code{PComponent} and its subclasses).
\end{enumerate}

\subsubsection{User-Visible Improvements}
\begin{enumerate}
\item Improving Auto (expanding its scope to givens).
\item Supporting long variable names, and possibly expanding the role of
variables along the lines of~\cite{Velleman06a}.
\item Supporting named hypothesis, to allow easy reference.
\item Allow proof comments.
\item Adding syntax highlighting (color coding of proofs).
\item Adding more toolbar buttons.
\item Updating HTML help files to reflect software changes.
\end{enumerate}

\subsection{To Be Done in Proof Maker}

Then our remaining work on Proof Maker includes the following.
\begin{enumerate}
\item Finishing and polishing the APM user interface (UI) as a genuine,
fully-functioning Android portal of the Proof Designer UI that has
the look-and-feel but also the behavior and user experience of native
Android applications.
\item Adding more touch-aware interactions to Proof Maker (e.g., dragging-and-dropping
of hypothesis, context menus).
\item Internationalizing Proof Maker, so as to allow languages such as Arabic,
Chinese, etc., in its proofs and its GUI.
\item Porting Proof Designer to other handheld device platforms such as
Windows 8 Phone and iOS.
\end{enumerate}

\section{Related Work}

Coq~\cite{Bertot2004} and Isabelle~\cite{Isabelle2015} are generic
proof assistants. Both build on a large tradition of scientific research
in the area of proof automation, going back to LCF~\cite{LCFPlotkin77,ML-LCF78}
and even further. Compared to Proof Designer, Coq and Isabelle are
vastly much more powerful (they can help construct proofs in almost
any mathematical domain), but the two proof assistants are much less
user-friendly than Proof Designer. Users of Coq and Isabelle have
to write code to construct their proofs. As such, Coq and Isabelle
users actually need to also be programmers, not only mathematicians
or math students. Since handheld devices typically lack a keyboard,
writing capabilities are usually limited on them. This casts doubts
on the likelihood of Coq or Isabelle getting ported to handheld devices.

DC Proof~\cite{DCProof2015} is a more user friendly software when
compared to Coq and Isabelle, and it is a bit more powerful than Proof
Designer. Given its ASCII-based mathematical notation, however, DC
Proof is less user-friendly than Proof Designer.

The following table summarizes some of the differences between these
proof construction software tools.

\noindent \begin{center}
\begin{tabular}{|c|c|c|c|c|}
\cline{2-5} 
\multicolumn{1}{c|}{} & PM & PD & DC Proof & Coq/Isabelle\tabularnewline
\hline 
Writing Code & No & No & No & Yes\tabularnewline
\hline 
\multirow{2}{*}{User-Friendly} & Yes (Touch- & \multirow{2}{*}{Yes} & Largely (Special & \multirow{2}{*}{No}\tabularnewline
 & based) &  & Math. Notation) & \tabularnewline
\hline 
Expressive & Elementary & Elementary & Elem. Set Theory, & Almost Any \tabularnewline
Power & Set Theory & Set Theory & Number Theory & Math Field\tabularnewline
\hline 
Automation & Very Limited & No & No & Yes\tabularnewline
\hline 
Saving Proofs & Yes & Export HTML & Yes & Yes\tabularnewline
\hline 
\end{tabular}
\par\end{center}

\section{Future Work}

Like Proof Designer, the scope of proofs doable in Proof Maker will
be limited to elementary set theory. Once done with Proof Maker, future
work that could be built on top of it can include adding a type system
(and possibly later a type inference system) that enables Proof Maker
overcome this fundamental limitation that it inherited from Proof
Designer. Adding a type system to Proof Maker will enable it to assist
in constructing proofs in mathematical domains other than set theory
(e.g., number theory, group theory, order theory, domain theory, etc.),
while maintaining the characteristic simplicity of the software and
its user-friendliness. (This software, for example, may help in our
formalization of an introductory domain theory textbook~\cite{AbdelGawad2015d}.)

Lurch is \textquotedblleft a word processor that can check your math\textquotedblright ~\cite{lurch}.
In particular, just as a word processor checks spelling and grammar
in natural language documents, Lurch aims to check any mathematical
proofs included in a document (e.g., a school math homework, an exam,
a research article, a book chapter, ... etc.) with as little user
guidance as possible (in the form of document annotations). In its
aims, Lurch was greatly influenced by Proof Designer. Adding a \emph{customizable}
type system to Lurch is another possible future work that can be done
after adding a type system to Proof Maker. Adding a customizable type
system to Lurch will add to Lurch the ability to restrict its rule
identifiers so that they can only be instantiated with an expression
of a certain type (like a statement, or a set, or a natural number,
... etc.). The type of an expression in Lurch will also be customizable,
and will be compatible with a customizable parser that the authors
of Lurch intend to soon add to Lurch.

\bibliographystyle{plain}

\newpage{}

\appendix
\noindent \begin{center}
\begin{figure}[H]
\noindent \begin{centering}
\includegraphics[bb=0bp 250bp 534bp 760bp,clip,scale=0.7]{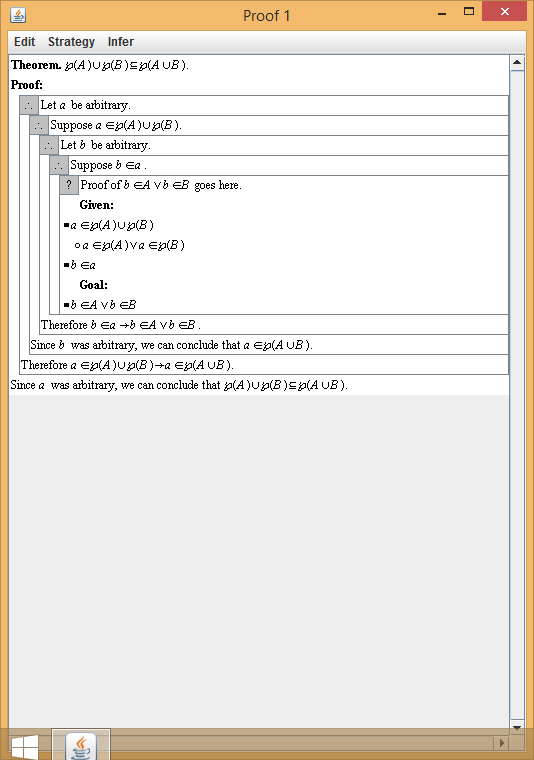}
\par\end{centering}

\protect\caption{\label{fig:An-Incomplete-Proof}An Incomplete Proof Designer Proof}
\end{figure}

\par\end{center}

\noindent \begin{center}
\begin{figure}[H]
\noindent \begin{centering}
\includegraphics[bb=2bp 0bp 502bp 558bp,clip,scale=0.7]{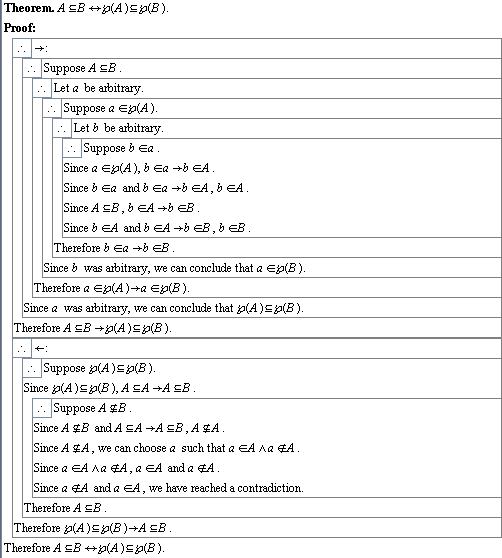}
\par\end{centering}

\protect\caption{\label{fig:A-Finished-Proof}A Finished Proof Designer Proof}
\end{figure}

\par\end{center}

\noindent \begin{center}
\begin{figure}[H]
\noindent \begin{centering}
\includegraphics[scale=0.7]{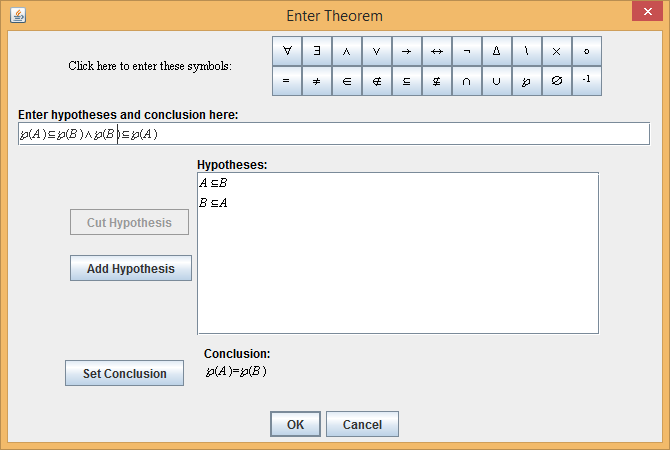}
\par\end{centering}

\protect\caption{\label{fig:Proof-Designer-Entry}Proof Designer Theorem Entry Dialog}
\end{figure}

\par\end{center}

\noindent \begin{center}
\begin{figure}[H]
\noindent \begin{centering}
\includegraphics[scale=0.7]{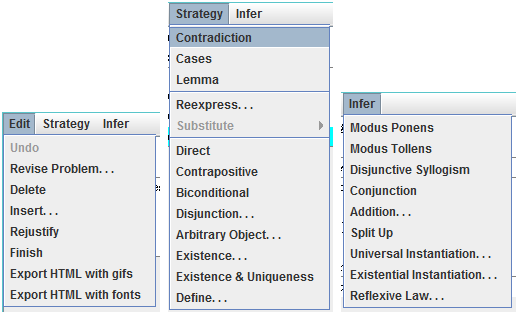}
\par\end{centering}

\protect\caption{\label{fig:Proof-Designer-Menus}Proof Designer Menus}
\end{figure}

\par\end{center}

\noindent \begin{center}
\begin{figure}[H]
\noindent \centering{}\includegraphics[scale=0.7]{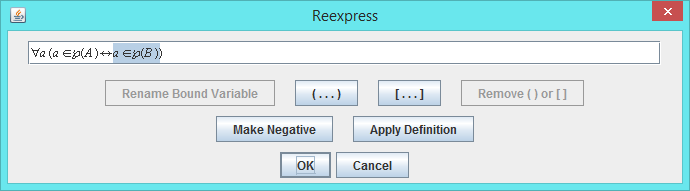}\protect\caption{\label{fig:Proof-Designer-Reexpress}Proof Designer Reexpress Dialog}
\end{figure}

\par\end{center}

\noindent \begin{center}
\begin{figure}[H]
\noindent \begin{centering}
\includegraphics[bb=200bp 550bp 900bp 700bp,clip,scale=0.5]{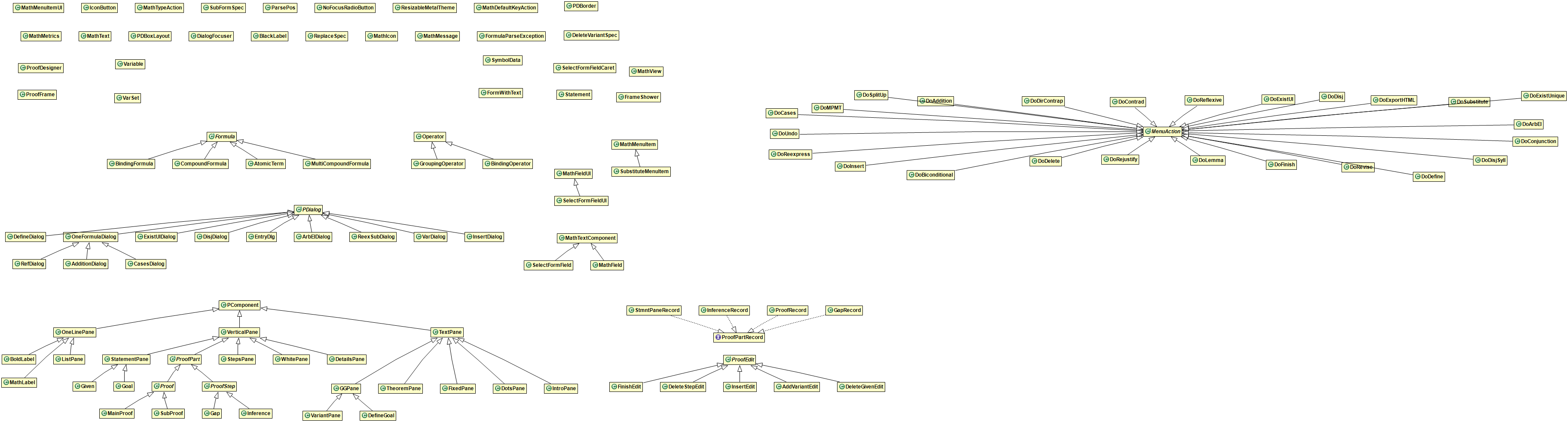}
\par\end{centering}

\noindent \begin{centering}
\includegraphics[bb=0bp 0bp 1350bp 350bp,clip,scale=0.3]{Class_Diagram__Inheritance_.png}
\par\end{centering}

\protect\caption{\label{fig:Class-hierarchy-PD}Class hierarchy for formulas and proof
components in Proof Designer}
\end{figure}

\par\end{center}

\noindent \begin{center}
\begin{figure}[H]
\noindent \begin{centering}
\includegraphics[scale=0.7]{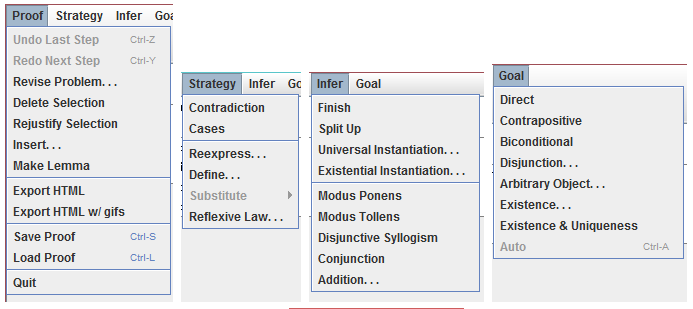}
\par\end{centering}

\protect\caption{\label{fig:NPD-Menus}New Proof Designer Menus}
\end{figure}

\par\end{center}

\noindent \begin{center}
\begin{figure}[H]
\noindent \begin{centering}
\includegraphics[scale=0.7]{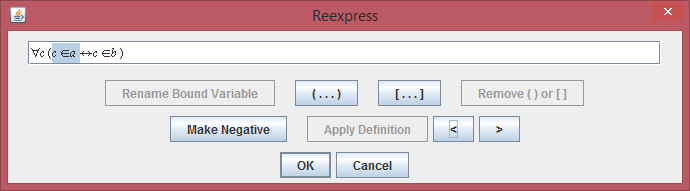}
\par\end{centering}

\protect\caption{\label{fig:NPD-Reexpress}New Proof Designer Reexpress Dialog (note
the < (undo) and > (redo) buttons)}
\end{figure}

\par\end{center}

\noindent \begin{center}
\begin{figure}[H]
\noindent \begin{centering}
\includegraphics[bb=0bp 335bp 534bp 537bp,clip,scale=0.7]{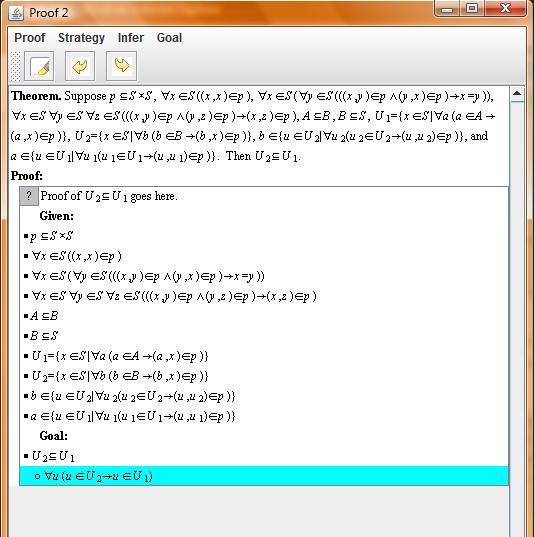}
\par\end{centering}

\protect\caption{\label{fig:NPD-Toolbar}New Proof Designer Toolbar}
\end{figure}

\par\end{center}

\noindent \begin{center}
\begin{figure}[H]
\noindent \begin{centering}
\includegraphics[bb=0bp 550bp 600bp 1024bp,clip,scale=0.5]{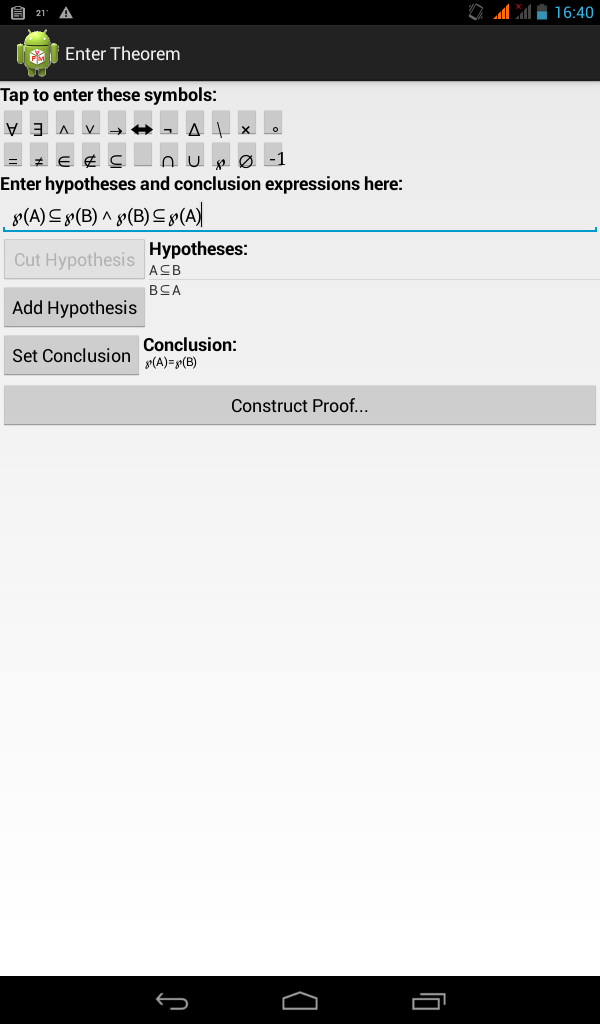}
\par\end{centering}

\protect\caption{\label{fig:APM-Theorem-Entry}APM Theorem Entry Dialog (proof of concept)}
\end{figure}

\par\end{center}

\noindent \begin{center}
\begin{figure}[H]
\noindent \begin{centering}
\includegraphics[bb=0bp 550bp 600bp 1024bp,clip,scale=0.5]{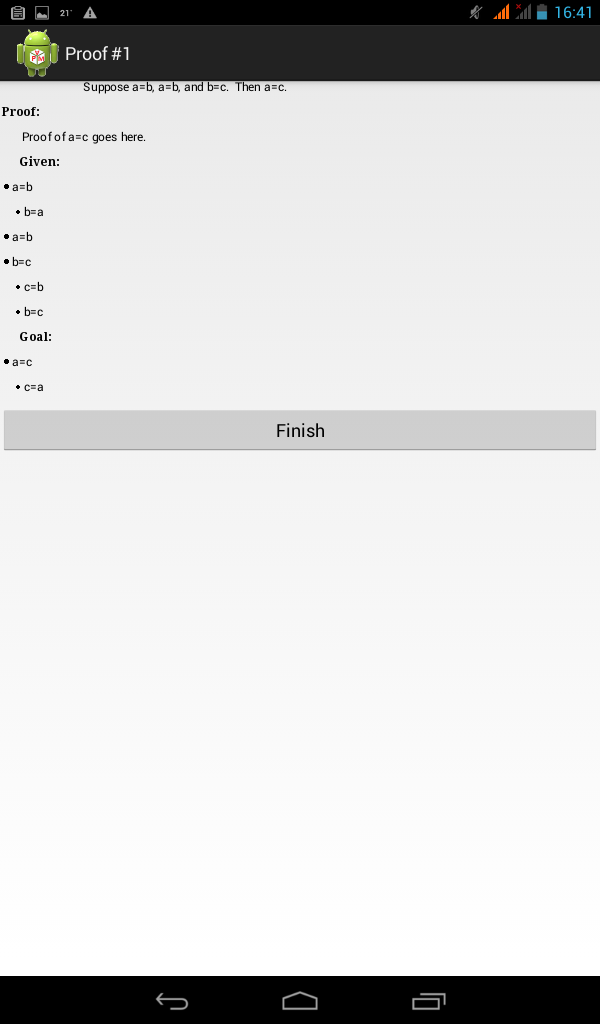}
\par\end{centering}

\protect\caption{\label{fig:APM-Incomplete-Proof}An Incomplete APM Proof (proof of
concept)}
\end{figure}

\par\end{center}
\end{document}